\input amstex
\documentstyle{amsppt}
\magnification=\magstep1
 \NoRunningHeads
 \NoBlackBoxes
\topmatter
\title
Krieger's type of nonsingular Poisson suspensions and IDPFT systems
\endtitle
\author Alexandre I. Danilenko and  Zemer Kosloff 
\endauthor

\abstract 
Given an infinite countable discrete amenable group $\Gamma$, we construct explicitly sharply weak mixing nonsingular Poisson $\Gamma$-actions of each Krieger's types: $III_\lambda$, for $\lambda\in[0,1]$, and $II_\infty$. 
The result is new even for $\Gamma=\Bbb Z$.
As these Poisson suspension actions are over very special dissipative base,
we obtain also new examples of
sharply weak mixing nonsingular Bernoulli $\Gamma$-actions and IDPFT systems
of each possible Krieger's type.
\endabstract

\address
B. Verkin Institute for Low Temperature Physics \& Engineering
of Ukrainian National Academy of Sciences,
47 Nauky Ave.,
 Kharkiv, 61164, UKRAINE
\endaddress
\email            alexandre.danilenko\@gmail.com
\endemail

\address
Einstein Institute of Mathematics, Hebrew University of Jerusalem, Givat Ram. Jerusalem, 9190401, ISRAEL
\endaddress
\email
zemer.kosloff\@mail.huji.ac.il
\endemail

\thanks
The research of Z.K. was partially supported by ISF grant No. 1570/17.
\endthanks

\loadbold

\endtopmatter

\document

\head 0. Introduction\endhead
Throughout this paper, $\Gamma$ denotes an  infinite countable discrete amenable group. 
Given an action of $\Gamma$ on an infinite volume measure space, the Poisson suspension of this action is the associated $\Gamma$-action on the
corresponding Poisson point process, whose elements can be viewed as 
clouds of particles without mutual interactions; the suspension  moves every  cloud particle-wise  according to the original action. 
See~Definition~1.4 below. 
In \cite{DaKoRo1} and \cite{DaKoRo2},  E. Roy and the present
authors initiated the study of nonsingular Poisson suspensions
and proved, among other things, that a generic (in the Baire category sense) nonsingular Poisson suspension is ergodic of Krieger's type $III_1$.
Our purpose in this work is to construct  examples of ergodic Poisson suspensions of other types (which are on the meager side).

\proclaim{Theorem 0.1}   
There exist sharply weak mixing nonsingular Poisson suspensions  of  totally dissipative  free $\Gamma$-actions of each of the following Krieger's types:  $II_\infty$ and $III_\lambda$, for  $\lambda\in[0,1]$.
\endproclaim

This result is new even in the case where $\Gamma=\Bbb Z$.
For the definition of sharp weak mixing (which is stronger than ergodicity and weak mixing) we refer to \S 1.

 It was shown  recently in \cite{KoSo} that for 
 each $\lambda\in(0,1)$, there exists a type $III_\lambda$ nonsingular Bernoulli  
 $\Gamma$-action.
 Bernoulli $\Gamma$-actions of type $II_\infty$ and $III_0$ were constructed in a later work
\cite{BeVa}.
Type $III_1$ Bernoulli actions of amenable groups appeared earlier in \cite{VaWa}.
We will show that each of the examples of Poisson suspensions in Theorem~0.1 is isomorphic to a nonsingular Bernoulli action of $\Gamma$.
Hence Theorem~0.1 gives a new  proof to    \cite{KoSo, Theorem~3}, 
 \cite{BeVa, Theorem~A} and the first claim of \cite{VaWa, Theorem~6.1}.

\proclaim{Corollary 0.2} For each  infinite countable discrete amenable group $\Gamma$, 
there exist sharply weak mixing nonsingular Bernoulli $\Gamma$-action
of each of the following Krieger's type:  $II_\infty$ and $III_\lambda$, for $\lambda\in[0,1]$.
\endproclaim

We also state one more byproduct of the proof of Theorem~0.1. 
It is about the Krieger's type for nonsingular  IDPFT (infinite direct product of finite types) systems, see Definition~1.1 below. 
 Such systems were introduced recently in \cite{DaLe}.
We note that only type $II_1$ and $III_1$ examples of  IDPFT systems have been known so far (see \cite{DaLe} and references therein).
We recall that D.~Ornstein and
B.~Weiss \cite{OrWe} proved that there exists
a unique (up to measure preserving isomorphism) infinite entropy probability preserving Bernoulli action of $\Gamma$.
This  is the $\Gamma$-shiftwise action  on the space $([0, 1], \text{Lebesgue})^\Gamma$.
 We will
refer to this action as the {\it infinite measure preserving Bernoulli
$\Gamma$-action}.

\proclaim{Corollary 0.3} Let $T=(T_\gamma)_{\gamma\in\Gamma}$ denote the measure preserving Bernoulli  $\Gamma$-action  of infinite entropy on a standard probability space $(X,\mu)$.
Then for each $\tau\in\{II_\infty\}\cup \{III_\lambda\mid 0\le\lambda\le 1\}$, there is a sequence
of probability measures $(\kappa_n)_{n=1}^\infty$ on $X$ such that $\kappa_n\sim\mu$ for each $n\in\Bbb N$,  the product $\bigotimes_{n=1}^\infty\kappa_n$ is quasi-invariant under $T^{\otimes \Bbb N}:=((T_\gamma)^{\otimes\Bbb N})_{\gamma\in\Gamma}$ and
the system $(X^\Bbb N, \bigotimes_{n=1}^\infty\kappa_n,T^{\otimes \Bbb N})$ is sharply weak mixing of Krieger type  $\tau$.
\endproclaim

The strategy of the proof of Theorem~0.1
 is as follows.
 We construct a certain family of nonsingular totally dissipative actions depending on parameters.
 A  quantitative criterion for conservativeness of nonsingular Poisson suspensions is provided 
 in Lemmata~1.5 and 1.6.
 Then it is shown that the Poisson suspensions of the totally dissipative actions under consideration have the structure of IDPFT  system. 
 Applying the aforementioned criterion we obtain that these Poisson suspensions are conservative.
 The conservativeness implies
 sharp weak mixing in view  of the IDPFT structure (see Corollary~1.7).
 This is done in Section~1.
It remains to compute the Krieger's type of these systems. 
We show that the type depends on concrete choice of parameters of the construction.
In Sections~2, 3, 4 and 5 we specify the parameters to obtain sharply weak mixing Poisson suspensions of Krieger's type $II_\infty$, $III_0$, $III_\lambda$ with $0<\lambda<1$ and $III_1$ respectively.

 {\sl Acknowledgement.}  We thank the anonymous referees for careful reading of the paper and helpful comments.

\head 1. Preliminaries and general construction\endhead

For a more detailed exposition of the following concepts (and proofs) we  refer the reader
to \cite{DaSi}, \cite{Aa}, \cite{Sc} and \cite{DaLe}.

\subhead{Weak topology}\endsubhead
 Let $(X,\goth B,\mu)$ be a $\sigma$-finite non-atomic standard measure space.
Denote by Aut$(X,\mu)$ the group of all $\mu$-nonsingular invertible transformations of $X$.
Given $S\in\text{Aut}(X,\mu)$, we denote by $U_S$ the unitary Koopman operator associated with $S$:
$$
U_Sf:=f\circ S^{-1}\sqrt{\frac{d\mu\circ S^{-1}}{d\mu}}\in\Cal U(L^2(X,\mu)),\qquad\text{for each }f\in L^2(X,\mu).
$$
The {\it weak} topology on Aut$(X,\mu)$ is induced by the weak operator topology on the unitary group
$\Cal U(L^2(X,\mu))$ via the embedding $S\mapsto U_S$.

\subhead{Conservativeness and weak mixing concepts}\endsubhead
A nonsingular action $T=(T(\gamma))_{\gamma\in\Gamma}$ of $\Gamma$ on $X$ (i.e. a group homomorphism $\Gamma\ni\gamma\mapsto T(\gamma)\in \text{Aut}(X,\mu)$) is called {\it conservative} if for each subset $B\in\goth B$ of positive measure, for a.e. $x\in B$, there is $\gamma\in\Gamma$ such that 
$T(\gamma )x\in B$.
We will use the following criterion of conservativeness: if $\mu(X)<\infty$ then $T$ is conservative if and only if $\sum_{\gamma\in\Gamma}\frac{d\mu\circ T(\gamma)}{d\mu}(x)=+\infty$ at a.e. $x\in X$.
If $T$ is free and there is a subset $B\in\goth B$ such that $X=\bigsqcup_{\gamma\in \Gamma}T(\gamma) B$ (mod 0)
then $T$ is called {\it totally dissipative}.
$T$ is called {\it ergodic} if each $T$-invariant Borel subset is either $\mu$-null or $\mu$-conull. 
If for each ergodic probability preserving $\Gamma$-action  $(S(\gamma))_{\gamma\in\Gamma}$, the direct product $(T(\gamma)\times S(\gamma))_{\gamma\in\Gamma}$ is ergodic then $T$ is called {\it weakly mixing}.
If $T$ is ergodic and for each   ergodic conservative nonsingular  $\Gamma$-action  $(S(\gamma))_{\gamma\in\Gamma}$, the direct product $(T(\gamma)\times S(\gamma))_{\gamma\in\Gamma}$ is ergodic whenever it is conservative then $T$ is called {\it sharply weak mixing}.

\subhead Krieger's type\endsubhead
Let $T$  be an ergodic $\Gamma$-action on $(X,\goth B,\mu)$.
Denote by $[T]$ the full group of $T$.
We recall that a transformation $R\in\text{Aut}(X,\mu)$ belongs to $[T]$
if and only if $Rx\in \{T(\gamma) x\mid \gamma\in \Gamma\}$ for a.e. $x\in X$.

If there is a $\mu$-equivalent $\sigma$-finite $T$-invariant measure then
$T$ is called {\it of type $II$}.
If the $T$-invariant measure is finite then $T$ is called {\it of type $II_1$};
if the $T$-invariant measure is infinite then $T$ is called {\it of type $II_\infty$}.
If $T$ is not of type $II$ then it is called {\it of type $III$}.
The type $III$ admits further classification into subtypes.
We first recall that
an element $r$ of the multiplicative group $\Bbb R^*_+$ is called an  {\it essential value of the Radon-Nikodym cocycle of $T$} if for each neighborhood $U$ of $r$ and each
subset $A\in\goth B$ of positive measure there exist a subset $B\in \goth B$ of positive measure and an element $\gamma\in\Gamma$ such that $B\cup T(\gamma) B\subset A$ and
$\frac{d\mu\circ T( \gamma)}{d\mu}(x)\in U$ for each $x\in B$.
The set of all essential values of the Radon-Nikodym cocycle of $T$ is denoted by $r(T)$.
It is a closed subgroup of $\Bbb R_+^*$.

If $r(T)=\Bbb R_+^*$ then $T$ is called of type $III_1$; if there is $\lambda\in(0,1)$ such that $r(T)=\{\lambda^n\mid n\in\Bbb Z\}$ then 
$T$ is called of type $III_\lambda$.
If $T$ is of type $III$ but not of type $III_\lambda$ for any $\lambda\in(0,1]$ then $T$ is called of type $III_0$.

We will also need the following folklore approximation result (see, for instance, \cite{ChHaPr, Lemma~2.2}).
Let $\goth B_0\subset \goth B$ be  a dense subalgebra. 
Let $\delta>0$ and $s\in\Bbb R_+^*$.
If for each $A\in\goth B_0$ of positive measure  and every neighborhood $U$ of $s$ there is a subset
$B\in \goth B$  and an element $R\in[T]$ such that $B\cup R B\subset A$, $\mu(B)>\delta\mu(A)$ and 
$\frac{d\mu\circ R}{d\mu}(x)\in U$ for each $x\in B$ then $s\in r(T)$.

\subhead Maharam extension and the associated flow\endsubhead
Let $\kappa$ denote the absolutely continuous measure on $\Bbb R$ such that
$d\kappa(t)=e^{-t}dt$.
Consider the product space $\widetilde X:=(X\times\Bbb R,\mu\otimes\kappa)$.
Given $\gamma\in\Gamma$ and $s\in\Bbb R$, define two transformations $\widetilde T(\gamma), \,R(s)\in\text{Aut}(\widetilde X,\mu\otimes\kappa)$ by setting for each $(x,t)\in\widetilde X$,
$$
\widetilde T(\gamma)(x,t):=\bigg(T(\gamma) x, t+\log\frac{d\mu\circ T(\gamma)}{d\mu}(x)\bigg)\quad\text{and}\quad
R(s)(x,t):=(x,t-s).
$$
Then $\widetilde T:=(\widetilde T(\gamma))_{\gamma\in\Gamma}$ is a $(\mu\otimes\kappa)$-preserving action of $\Gamma$ on $\widetilde X$  
and 
$R:=(R(s))_{s\in\Bbb R}$ is a $(\mu\otimes\kappa)$-nonsingular action  of $\Bbb R$ on $\widetilde X$.
Moreover, $\widetilde T(\gamma)R(s)=R(s)\widetilde T(\gamma)$
for all $\gamma\in\Gamma$ and $s\in\Bbb R$.
Hence the restriction of  
$R$ to the $\sigma$-algebra of 
$\widetilde T$-invariant subsets  of $X\times\Bbb R$ equipped with (the restriction of) $\mu\otimes\kappa$ is well defined as a nonsingular action.
It is always ergodic (under assumption  that $T$ is ergodic).
It is called {\it the associated flow} of $T$ and denoted by $W^T$.
The action $\widetilde T$ is called {\it the Maharam ($\mu$-skew product) extension} of $T$.
The following holds:
\roster
\item"---" $T$ is of type $II$ if and only if $W^T$ is  transitive and aperiodic.
\item"---" $T$ is of type $III_\lambda$ with $0<\lambda<1$ if $W^T$ is transitive but periodic with period $\log\lambda$ .
\item"---" $T$ is of type $III_1$ if $W^T$ is the trivial action on a singleton.
\item"---" $T$ is of type $III_0$ if $W^T$ is non-transitive.
\endroster

\subhead IDPFT systems\endsubhead
The following definition generalizes naturally the concept of IDPFT introduced in \cite{DaLe}  for $\Bbb Z$-actions to actions of arbitrary countable groups.

\definition{Definition 1.1} Let $S_n$ be a nonsingular $\Gamma$-action on a standard probability space $(Y_n,\goth C_n,\mu_n)$ for each $n\in\Bbb N$ and let
$(Y,\goth C,\mu):=\bigotimes_{n\in\Bbb N}(Y_n,\goth C_n,\mu_n)$.
Suppose that the infinite direct product $S:=\bigotimes_{n\in\Bbb N} S_n$ is a $\mu$-nonsingular $\Gamma$-action.
If for each $n\in\Bbb N$, there exists an $S_n$-invariant $\mu_n$-equivalent probability measure $\nu_n$ on $Y_n$ then the dynamical system 
$(Y,\goth C,\mu, S)$ is called {\it an infinite direct product of finite types (IDPFT)}.
Moreover, if we set $f_n:=\frac{d\mu_n}{d\nu_n}$ then for each $\gamma\in\Gamma$,
$$
\frac{d\mu\circ S(\gamma)}{d\mu}(y)=\prod_{n=1}^\infty\frac{f_n(S_n(\gamma)y_n)}{f_n(y_n)}\qquad\text{at $\mu$-a.e. $y=(y_n)_{n=1}^\infty\in \prod_{n=1}^\infty Y_n$.}\tag1-1
$$
\enddefinition

We state the following claim without proof since it is a minor modification of \cite{DaLe, Proposition~2.3} proved there in the case where $\Gamma$ is isomorphic to $\Bbb Z$: the point is that the Schmidt-Walters  theorem (see \cite{ScWa} and \cite{DaLe, Theorem~B}) holds in this generality.

\proclaim{Fact 1.2} Let $S$ be an IDPFT system as in Definition~1.1.
Let
$(Y_n,\goth C_n,\mu_n)$ be mildly mixing for each $n\in\Bbb N$.
If $S$ is $\mu$-conservative then $S$ is 
$\mu$-sharply weak mixing.
\endproclaim

Denote by $\Psi$ the action of the group $\bigoplus_{n=1}^\infty \Gamma$ on  $(Y,\mu)$ generated by the transformations $\Psi(\gamma_1,\dots,\gamma_n):=S_1(\gamma_1)\times\cdots\times S_n(\gamma_n)\times I\in \text{Aut}(Y,\mu)$, $\gamma_1,\dots,\gamma_n\in\Gamma$, $n\in\Bbb N$.
Then the Maharam extension $\widetilde\Psi$ of $\Psi$ is defined
on the same space as the Maharam extension $\widetilde S$ of $S$.

\proclaim{Proposition 1.3}  Let $S$ be an IDPFT system as in Definition~1.1.
Let
$(Y_n,\goth C_n,\mu_n)$ be mildly mixing for each $n\in\Bbb N$.
If $(Y,\mu, S)$  is conservative 
then the $\sigma$-algebra $\goth I(\widetilde {S})$ of $\widetilde {S}$-invariant measurable subsets coincides with the $\sigma$-algebra $\goth I(\widetilde \Psi)$ of $\widetilde \Psi$-invariant subsets.
Hence the flow associated with $S$ coincides with the flow associated with $\Psi$.
\endproclaim

\demo{Proof}
We first claim that
 $\goth I(\widetilde{ S})\subset \goth I(\widetilde \Psi)$  (cf. Claim~I of \cite{DaLe, Theorem~2.10} in the case where $\Gamma$ is isomorphic to $\Bbb Z$).
 Take a subset $A\in \goth I(\widetilde{ S})$.
 Fix $n\in\Bbb N$.
 The measure $\mu^{(n)}:=\big(\bigotimes_{k=1}^n\nu_k\big)\otimes\bigotimes_{k>n}\mu_k$
 is equivalent to $\mu$.
 Denote by $\widehat S$ and $\widehat\Psi$ the Maharam extension of $S$ and $\Psi$ respectively with respect to $\mu^{(n)}$.
 Denote by $\xi$ the mapping
 $$
 \xi:X\times\Bbb R\ni (x,t)\mapsto \xi(x,t):=\Big(x, t+\log\frac{d\mu^{(n)}}{d\mu}(x)\Big)\in X\times\Bbb R.
 $$
 It is straightforward to verify that  $\xi$ in an isomorphism of $ (Y\times\Bbb R,\mu^{(n)}\otimes\tau)$ onto 
 $ (Y\times\Bbb R,\mu\otimes\tau)$ such that $\xi^{-1}\widetilde S(\gamma)\xi=\widehat S(\gamma)$ for each $\gamma\in\Gamma$ and
 $\xi^{-1}\widetilde\Psi(\gamma_1,\dots,\gamma_n)\xi=\widehat \Psi(\gamma_1,\dots,\gamma_n)$ for all $\gamma_1,\dots,\gamma_n\in\Gamma$.
Moreover,  there is  a nonsingular $\Gamma$-action $M=(M(\gamma))_{\gamma\in\Gamma}$ on the space $\big(\big(\bigotimes_{k>n}Y_k\big)\times\Bbb R, \big(\bigotimes_{k>n}\mu_k\big)\otimes\tau\big)$ such that
$\widehat S$ splits into the direct product in the following way:
$\widehat S(\gamma)=\big(\bigotimes_{k=1}^n S_k(\gamma)\big)\otimes M(\gamma)$ for each $\gamma\in\Gamma$.
Since $A$ is invariant under $\widetilde S$, it follows that the subset $\xi^{-1} A\subset Y\times\Bbb R$ is invariant under $\widehat S$.
Since $\widehat S$ is conservative (as $\widetilde S$ is conservative), and the $\Gamma$-action
$\bigotimes_{k=1}^n S_k$ on  the space $\bigotimes_{k=1}^n(Y_k,\mu_k\big)$  is mildly mixing, it follows from \cite{ScWa} that there is a subset $B\subset \big(\bigotimes_{k>n}Y_k\big)\times\Bbb R$ such that $\xi^{-1} A=\big(\bigotimes_{k=1}^n Y_k\big)\times B$.
It follows that  $\xi^{-1} A$ is invariant under $\big(\bigotimes_{k=1}^n S_k(\gamma_k)\big)\otimes I=\widehat\Psi(\gamma_1,\dots,\gamma_k)$ for all $\gamma_1,\dots,\gamma_n\in\Gamma$.
Hence $A$ is invariant under $\widetilde\Psi(\gamma_1,\dots,\gamma_k)$ for all $\gamma_1,\dots,\gamma_n\in\Gamma$. 
Since $n$ is arbitrary, the claim is proved. 

We now prove the opposite inclusion: $ \goth I(\widetilde \Psi)\subset\goth I(\widetilde{ S})$.
We note that  for each $\gamma\in\Gamma$,
$$
S(\gamma)=\lim_{n\to\infty}\Psi(\,\underbrace{\gamma,\dots,\gamma}_{n}\,)\quad
 \text{in the weak topology.}
 $$
It follows that $
\widetilde S(\gamma)=\lim_{n\to\infty}\widetilde \Psi(\,\underbrace{\gamma,\dots,\gamma}_{n}\,)$ in the weak topology.
Hence the transformation $\widetilde S(\gamma)$ is contained in the weak closure of $\{\widetilde\Psi(\theta)\mid\theta\in \bigoplus_{n=1}^\infty\Gamma\}$.
The desired inclusion  follows.
\qed
\enddemo

\subhead Nonsingular Poisson suspensions\endsubhead
 Let $(X,\goth B)$ be a standard Borel space and let $\mu$ be an infinite $\sigma$-finite 
non-atomic measure on $X$.
 Let $X^*$ be the set of purely atomic  ($\sigma$-finite) measures on $X$.
  For each subset $A\in\goth B$ with $0<\mu(A)<\infty$, we define a mapping $N_A:X^*\to\Bbb R$ by setting
 $N_A(\omega):=\omega(A)$.
 Let $\goth B^*$ stand for the smallest $\sigma$-algebra on $X^*$ such that the mappings
 $N_A$ are all $\goth B^*$-measurable.
There is a unique probability measure $\mu^*$  on $(X^*,\goth B^*)$ satisfying the following two conditions:
\roster
\item"---"  the measure $\mu^*\circ N_A^{-1}$ is the Poisson distribution
with parameter $\mu(A)$ for each $A\in\goth B$ with $0<\mu(A)<\infty$,
\item"---" given a countable family $A_1,A_2,\dots$ of mutually disjoint subsets $A_1,A_2,\dots\in \goth B$ with $0<\mu(A_n)<\infty$ for each $n$, the corresponding random variables $N_{A_1},N_{A_2},\dots, $ defined on the space $(X^*,\goth B^*, \mu^*)$ are independent.
\endroster
 Then $(X^*,\goth B^*, \mu^*)$ is a probability Lebesgue space.
 It is called {\it the Poisson point process on $X$ with intensity measure $\mu$}. 
Since $\mu$ is non-atomic,
$\mu^*$
is a {\it simple} point process, i.e.  for $\mu^*$-almost every
$\omega\in X^*$, we have that $\omega(\{x\})\in\{0,1\}$
 for all $x\in X$.
We let 
$$
\text{Aut}_1(X,\mu):=\bigg\{R\in \text{Aut}(X,\mu)\mid {\frac{d\mu\circ R}{d\mu}}-1\in L^1(X,\mu)\bigg\}.
$$ 
If $R\in \text{Aut}_1(X,\mu)$, we put $\chi(R):=\int_X( {\frac{d\mu\circ R}{d\mu}}-1)d\mu$.
Then Aut$_1(X,\mu)$ is  a subgroup of Aut$(X,\mu)$ and $\chi$ is a homomorphism
of Aut$_1(X,\mu)$ onto $\Bbb R$.
Given a transformation $R\in \text{Aut}_1(X,\mu)$, we define a map $R_*:X^*\to X^*$ by setting $R_*\omega:=\omega\circ R^{-1}$ for each $\omega\in X^*$.
Then  $R_*\in\text{Aut}(X^*,\mu^*)$ \cite{DaKoRo1, \S 4} and for $\mu^*$-a.e. $\omega$,
$$
\frac{d\mu^*\circ R_*}{d\mu^*}(\omega)=e^{-\chi(T)}\prod_{\omega(\{x\})>0}\frac{d\mu\circ R}{d\mu}(x).
\tag1-2
$$
\definition{Definition 1.4} The transformation $R_*$ is called  {\it the (nonsingular) Poisson suspension} of $R$.
More generally, given an action $T=(T(\gamma))_{\gamma\in\Gamma}$ of $\Gamma$ such that $T(\gamma)\in\text{Aut}_1(X,\mu)$,
we call the $\Gamma$-action $T_*=(T(\gamma)_*)_{\gamma\in\Gamma}$ 
{\it the (nonsingular) Poisson suspension} of $T$.
\enddefinition

We now show how to deduce Corollary~0.2 from Theorem~0.1.
If $T=(T_\gamma)_{\gamma\in\Gamma}$ is a free totally dissipative  $\Gamma$-action on
a measure space $(X,\mu)$ then there is a subset $A\subset X$ of positive measure such that $X=\bigsqcup_{\gamma\in\Gamma} T_\gamma A$.
Suppose that the nonsingular  Poisson suspension $(X^*,\mu^*, T_*)$ of $(X,\mu, T)$ is well defined.
Then by the  property of independence for the Poisson suspensions, the system $(X^*,\mu^*, T_*)$ is canonically isomorphic to the left shift-wise action of $\Gamma$ on the infinite product measure space
$\prod_{\gamma\in G}(A^*,(\mu\circ T_\gamma^{-1})^*)$.
In other words, $(X^*,\mu^*, T_*)$ is a nonsingular Bernoulli $\Gamma$-action.
Therefore Corollary~0.2  follows from Theorem~0.1.

Given  a subset $B\in\goth B$ and an integer $n\in\Bbb Z_+$, we denote by
$[B]_n$ the cylinder $\{\omega\in X^*\mid \omega(B)=n\}$.
If $\nu$ is a $\sigma$-finite measure on $(X,\goth B)$ then $\nu^*\sim\mu^*$ if and only if
$\mu\sim\nu$ and $\sqrt{\frac{d\mu}{d\nu}}-1\in L^2(X,\nu)$ (see \cite{Ta} or \cite{DaKoRo1, Theorem~3.3}).
We will use below the  following sufficient condition for conservativeness, which extends
  \cite{DaKoRo2, Proposition~4.4} to arbitrary countable group actions.

\proclaim{Lemma 1.5} Let \,$T=(T(\gamma))_{\gamma\in\Gamma}$ be an action of $\Gamma$ such that  $T(\gamma)\in\text{\rom{Aut}}_1(X,\mu)$, $\chi(T(\gamma))=0$  and
$\Big(\frac {d\mu}{d\mu\circ T(\gamma)}\Big)^2-1\in L^1(X,\mu)$  for each $\gamma\in\Gamma$.
If there is a sequence $(b_\gamma)_{\gamma\in\Gamma}$ of positive  reals such that
$\sum_{\gamma\in\Gamma} b_\gamma=\infty$ but $\sum_{\gamma\in\Gamma} 
b_\gamma^2 
e^{\int_X \big(\big(\frac {d\mu}{d\mu\circ T(\gamma)}\big)^2-1\big) d\mu}<\infty$ then 
$T_*$ is conservative.
\endproclaim

\demo{Proof}
As in the proof of \cite{DaKoRo2, Proposition~4.4}, it follows from the 
 assumptions of the lemma that for each $\gamma\in\Gamma$,
$$
M_\gamma :=
\bigg\|\frac{d\mu^*}{d\mu^*\circ T(\gamma)_*}\bigg\|_2^2 =
e^{\int_X\big(\big( \frac{d\mu}{d\mu\circ T(\gamma)}\big)^2-1\big)d\mu}.
$$
By Markov's inequality, 
$$
\mu^*\bigg(\bigg\{\omega\in X^*\mid \frac{d\mu^*}{d\mu^*\circ T(\gamma)_*}(\omega)>\frac1{b_\gamma}
\bigg\}\bigg)\le b_\gamma^2 M_\gamma=b_\gamma^2 e^{\int_X\big(\big( \frac{d\mu}{d\mu\circ T(\gamma)}\big)^2-1\big)d\mu}.
$$
As the righthand side is summable, it follows from the Borel-Cantelli lemma  that $
\frac{d\mu^*}{d\mu^*\circ T(\gamma)_*}(\omega)\le\frac1{b_\gamma}$ for all but finitely many $\gamma\in\Gamma$ at a.e. $\omega$.
The latter inequality is equivalent to $\frac{d\mu^*\circ T(\gamma)_*}{d\mu^*}(\omega)\ge b_\gamma$.
Since $\sum_{\gamma\in\Gamma}b_\gamma=\infty$, we conclude  that $\sum_{\gamma\in\Gamma}\frac{d\mu^*\circ T( \gamma)_*}{d\mu^*}(\omega)=+\infty$ at a.e. $\omega$.
Hence $T_*$ is conservative.
\qed
\enddemo

\subhead{General construction}\endsubhead
Let $\Gamma$ be a discrete  infinite  countable group.
Fix two sequences $(a_n)_{n=1}^\infty$ and $(\lambda_n)_{n=1}^\infty$ of positive reals and a sequence $(F_n)_{n=1}^\infty$
of finite subsets in $\Gamma$.
Denote by $m_\Gamma$ the counting measure on $\Gamma$.
Denote by $m_\Bbb T$ the Haar measure on the circle $\Bbb T$.
Let $Z:=\Bbb T\times\Gamma$.
We define two mutually commuting Borel $\Gamma$-actions $L=(L(\gamma))_{\gamma\in\Gamma}$ and $R=(R(\gamma))_{\gamma\in\Gamma}$ 
by setting
$$
L(\gamma)(t,r):=(t,\gamma r),\quad R(\gamma)(t,r):=(t, r\gamma^{-1})
\qquad\text{for each }(t,r)\in Z.\tag1-3
$$
Define  two sequences of measures $(\nu_n)_{n=1}^\infty$ and $(\mu_n)_{n=1}^\infty$ on $Z$ by setting
$\nu_n:=\frac{a_n}{\# F_n}m_\Bbb T\otimes m_\Gamma$ and $\mu_n\sim\nu_n$ with
$$
\frac{d\mu_n}{d\nu_n}(t,r):=1+(\lambda_n-1)1_{F_n}(r)\qquad\text{for each }(t,r)\in Z.
$$
Then  $L$ and $R$  are  both totally dissipative, free, $\mu_n$-nonsingular and $\nu_n$-preserving for each $n\in\Bbb N$.
We now let $(X,\nu,T):=\bigsqcup_{n=1}^\infty(Z,\nu_n,L)$ and
$\mu:=\bigsqcup_{n=1}^\infty\mu_n$.
Then the $\Gamma$-action $T$ preserves
$\nu$ and leaves $\mu$ quasi-invariant.
In particular, we can write that $T(\gamma)\in \text{Aut}_1(X,\nu)$ for each $\gamma\in\Gamma$.

For each $\lambda>0$ and $\lambda\ne 1$, we let
$c(\lambda):=\lambda^3-\lambda+\lambda^{-2}-1=(1-\lambda^2)(1-\lambda^3)\lambda^{-2}>0$.
The following lemma (jointly with Lemma~1.5) will be used  in subsequent sections to check conservativeness of the nonsingular Poisson
suspensions $T_*$.

\proclaim{Lemma 1.6} Fix $\gamma\in\Gamma$.
\roster
\item"\rom{(i)}"
$T(\gamma)\in\text{\rom{Aut}}_1(X,\mu)$ if and only if $\sum_{n=1}^\infty \frac{\#(\gamma F_n\triangle F_n)}{\# F_n}|1-\lambda_n|a_n<+\infty$.
Moreover, $\chi(T(\gamma))=0$ whenever  $T(\gamma)\in\text{\rom{Aut}}_1(X,\mu)$.
\item"\rom{(ii)}"
$ \big(\frac {d\mu}{d\mu\circ T(\gamma)^{-1}}\big)^2-1\in L^1(X,\mu)$ if and only if 
$$\sum_{n=1}^\infty 
\frac{|1-\lambda_n^2|(1+\lambda_n^3)}{\lambda_n^2}\frac{\#(\gamma F_n\triangle F_n)}{\# F_n}a_n<+\infty.
$$
Moreover,
$\int_X\Big(\big(\frac {d\mu}{d\mu\circ T(\gamma)^{-1}}\big)^2-1\Big) d\mu=\sum_{n=1}^{\infty} c(\lambda_n)a_n\frac{\#(\gamma F_n\triangle F_n)}{2\# F_n}$.
\endroster
\endproclaim
\demo{Proof}
(i) Let $f_n:=\frac{d\mu_n}{d\nu_n}$.
For each $\gamma\in\Gamma$ and $n\in\Bbb N$,
$$
\align
\int_X\bigg|\frac{d\mu\circ T(\gamma)^{-1}}{d\mu}&-1\bigg|d\mu
=\sum_{n=1}^\infty\int_Z \bigg|\frac{d\mu_n\circ L(\gamma)^{-1}}{d\mu_n}-1\bigg|    d\mu_n\\
&=\sum_{n=1}^\infty\int_Z|f_n\circ L(\gamma)^{-1}-f_n|d\nu_n
\\
&=\sum_{n=1}^\infty|1-\lambda_n|\frac{a_n}{\# F_n}\int_{F_n\cup\gamma F_n}|1_{F_n}\circ \gamma^{-1}-1_{F_n}|\,dm_\Gamma
\\
&=\sum_{n=1}^\infty\frac{\#(\gamma F_n\triangle F_n)}{\# F_n}|1-\lambda_n|a_n,
\endalign
$$
and the first claim of Lemma~1.6(i) is proved.
In a similar way, 
$$
\align
\chi(T(\gamma))&:=\int_X\bigg(\frac{d\mu\circ T(\gamma)^{-1}}{d\mu}-1\bigg)d\mu \\
&=\sum_{n=1}^\infty\int_{Z}(f_n\circ L(\gamma)^{-1}-f_n)d\nu_n\\
&
=\sum_{n=1}^\infty\frac{a_n(\lambda_n-1)}{\# F_n}\bigg(\int_{F_n\cup\gamma F_n}1_{F_n}\circ \gamma^{-1}dm_\Gamma-\int_{F_n\cup\gamma F_n}1_{F_n}dm_\Gamma\bigg)=
0.
\endalign
$$
Thus,~Lemma~1.6(i) is proved completely.

(ii) 
We
 first observe that
$$
\int_X\bigg|\bigg(\frac {d\mu}{d\mu\circ T(\gamma)^{-1}}\bigg)^2-1\bigg| d\mu
=\sum_{n=1}^\infty\int_{Z}\bigg|\frac{f_n^3}{(f_n\circ L(\gamma)^{-1})^2}-f_n\bigg|d\nu_n.
$$
Next, we see that
$$
\align
\int_{Z}\bigg|\frac{f_n^3}{(f_n\circ L( \gamma)^{-1})^2}-f_n\bigg|d\nu_n
&=|\lambda_n^3-\lambda_n|\frac{a_nm_\Gamma(F_n\setminus \gamma F_n)}{\# F_n}+|\lambda_n^{-2}-1|\frac{a_nm_\Gamma(\gamma F_n\setminus F_n)}{\# F_n}\\
&= \frac{|1-\lambda_n^2|(1+\lambda_n^3)}{\lambda_n^2}\frac{a_n\#(\gamma F_n\triangle F_n)}{2\# F_n}
\endalign
$$
and
the first claim of (ii) follows.
In a similar way,
$$
\align
\int_X\bigg(\bigg(\frac {d\mu}{d\mu\circ T(\gamma)^{-1}}\bigg)^2-1\bigg) d\mu
&=\sum_{n=1}^\infty\bigg((\lambda_n^3-\lambda_n)\frac{a_nm_\Gamma(F_n\setminus \gamma F_n)}{\# F_n}\\
&+(\lambda_n^{-2}-1)\frac{a_nm_\Gamma(\gamma F_n\setminus F_n)}{\# F_n}\bigg)\\
&=\sum_{n=1}^\infty c(\lambda_n)\frac{\#(\gamma F_n\triangle F_n)}{2\# F_n}a_n,
\endalign
$$
as desired.
\qed
\enddemo

\subhead IDPFT structure of the Poisson suspensions
\endsubhead
We now describe a special infinite product structure of the Poisson suspensions of the systems considered above.
Enumerate the elements of $\Gamma$: $\Gamma=(\gamma_n)_{n=1}^\infty$.
Select a F{\o}lner sequence $(F_n)_{n=1}^\infty$ in $\Gamma$ such that
$$
\max_{1\le k\le n}\frac{\#(\gamma_kF_n\triangle F_n)}{\# F_n}\le\frac 1{n^2}\qquad\text{for  each $n\in\Bbb N$.}
\tag1-4
$$
In all constructions below, the enumeration of $\Gamma$ and the F{\o}lner sequence are
as above. 
We recall that $(X,\mu,T)=\bigsqcup_{n=1}^\infty(Z,\mu_n,L)$.
Thus, $X$ is a union of countably many disjoint $T$-invariant enumerated ``copies'' of $Z$.
Given $\omega\in X^*$, consider the sequence of restrictions of $\omega$ to these copies.
We thus obtain a one-to-one mapping from $X^*$ to $(Z^*)^\Bbb N$.
The independence property of Poisson suspensions yields that this mapping maps $\mu^*$ onto the direct product $\bigotimes_{n\in\Bbb N}\mu_n^*$. 
Thus, we obtain the canonical isomorphism of $(X^*,\mu^*,T_*)$ onto $\bigotimes_{n=1}^\infty(Z^*,\mu_n^*,L_*)$.
If $T(\gamma)\in\text{Aut}_1(X,\mu)$ for each $\gamma\in\Gamma$, the Poisson suspension $T_*$ is $\mu^*$-nonsingular.
Hence the infinite product system $\bigotimes_{n=1}^\infty(Z^*,\mu_n^*,L_*)$ is also nonsingular.
Since $\mu_n\sim\nu_n$ and $\sqrt{\frac{d\mu_n}{d\nu_n}}-1\in L^2(Z,\nu_n)$, the dynamical system $(Z^*,\mu_n^*,L_*)$ is of Krieger's type $II_1$ for each
$n\in\Bbb N$.
The corresponding $L_*$-invariant $\mu_n^*$-equivalent probability measure is $\nu_n^*$.

Thus,  we showed that if $T(\gamma)\in\text{Aut}_1(X,\mu)$ for each $\gamma\in\Gamma$,
then $(X^*,\mu^*, \Gamma_*)$ is (isomorphic to) an IDPFT system.
Since $L$ is totally dissipative, free   and $\nu_n$-preserving, the system $(Z^*,\nu_n^*,L_*)$ is a measure preserving Bernoulli action of $\Gamma$ for each $n\in\Bbb N$.
It is mixing (hence mildly mixing).
Therefore we deduce the following corollary from Fact~1.2.

\proclaim{Corollary 1.7} If the Poisson suspension $(X^*,\mu^*, T_*)$ is nonsingular and conservative then
it  is sharply weak mixing.
\endproclaim

Let $A_n:=\Bbb T\times F_n$, $n\in\Bbb N$.
It follows from \cite{DaKoRo1, Theorem~3.6} that for each $n\in\Bbb N$ and $\mu_n^*$-a.e. $\omega\in Z^*$,
$$
\aligned
\frac{d\nu_n^*}{d\mu_n^*}(\omega)&=e^{-\int_{Z}(\frac 1{f_n}-1)d\mu_n}\prod_{\{x\in Z\mid\omega(\{x\})=1\}}\frac{d\nu_n}{d\mu_n}(x)\\
&=e^{\mu_n(A_n)-\nu_n(A_n)}\lambda_n^{-\omega(A_n)}\\
&=e^{\mu_n(A_n)-\nu_n(A_n)}\sum_{k=0}^\infty\lambda_n^{-k}1_{[A_n]_k}(\omega).
\endaligned
\tag1-5
$$

\head 2. Poisson suspensions of type $II_\infty$
\endhead

We first recall two concepts  of $\sigma$-finite  products for  a sequence of probability spaces (see \cite{Hi} and \cite{Mo}).

Throughout this section  $(Y_n,\goth C_n, \eta_n)$ is a standard probability space for all  $n\in\Bbb N$.

\definition{Definition 2.1} We call a  $\sigma$-finite measure $\eta$ on the standard Borel space $(Y,\goth C):=\bigotimes_{n\in\Bbb N}(Y_n,\goth C_n)$  a {\it Moore-Hill-product (MH-product)} of $(\eta_n)_{n=1}^\infty$ if for each $n\in\Bbb N$, there exists a $\sigma$-finite measure $\eta_n^*$ on the infinite product space 
$\bigotimes_{k>n}(Y_k,\goth C_k)$ such that $\eta=\eta_1\otimes\cdots\otimes\eta_n\otimes\eta_n^*$.
\enddefinition

For each $n\in\Bbb N$, fix a subset $B_n\in \goth C_n$ such that $\eta_n(B_n)>0$.
Denote the sequence $(B_n)_{n=1}^\infty$ by $\boldkey B$.
Let $\boldkey B^n:=Y_1\times\cdots\times Y_n\times B_{n+1}\times B_{n+2}\times\cdots\in\goth B$.
Then  $\boldkey B^1\subset \boldkey B^2\subset\cdots$.
A $\sigma$-finite measure $\eta^{\boldkey B}$ on $(Y,\goth C)$ is well defined by the following sequence of restrictions:
$$
\eta^{\boldkey B}\restriction\boldkey B^n:=\frac{\eta_1}{\eta_1(B_1)}\otimes\cdots\otimes\frac{\eta_n}{\eta_n(B_n)}
\otimes\frac{\eta_{n+1}\restriction B_{n+1}}{\eta_{n+1}(B_{n+1})} \otimes\frac{\eta_{n+1}\restriction B_{n+1}}{\eta_{n+1}(B_{n+1})}\otimes\cdots,
$$
for each $n\in\Bbb N$.

\definition{Definition 2.2 \cite{Hi, \S3.1}} $\eta^{\boldkey B}$ is called {\it the restricted product} of $(\eta_n)_{n=1}^\infty$ with respect to $\boldkey B$.
\enddefinition

We note that $\eta^\boldkey B$ is supported on  $\bigcup_{n=1}^\infty\boldkey B^n$ and $\eta^\boldkey B(\boldkey B^n)=\prod_{j=1}^n\eta_n(B_n)^{-1}$ for each $n\in\Bbb N$.
This  and \cite{Hi, Corollary~3.5} implies the following corollary.

\proclaim{Corollary 2.3} The following are equivalent:
\roster
\item"\rom{(i)}"
 $\eta^{\boldkey B}$ is finite,
\item"\rom{(ii)}"
 $\prod_{n=1}^\infty\eta_n(B_n)>0$,
  \item"\rom{(iii)}"
  $\eta^{\boldkey B}\sim\bigotimes_{n=1}^\infty\eta_n$.
  \endroster
 \endproclaim

Given two probability measures $\alpha,\beta$ on a standard Borel space $(Y,\goth C)$,
let $\eta$ be a third probability measure on $\goth C$ such that $\alpha\prec\eta$ and $\beta\prec\eta$.
The (squared) {\it Hellinger distance between $\alpha$ and $\beta$} is
$$
H^2(\alpha,\beta):=\frac12\int_Y\Bigg(\sqrt{\frac{d\alpha}{d\eta}}-\sqrt{\frac{d\beta}{d\eta}}\Bigg)^2
d\eta
$$
This definition does not depend on the choice of $\eta$.
We will utilize the following results from \cite{Hi}.

\proclaim{Fact 2.4} Let $\eta_n$ and $\mu_n$ be two probability measures on $(Y_n,\goth C_n)$ for each $n\in\Bbb N$.
Let $\boldkey B=(B_n)_{n=1}^\infty$ be a sequence of  subsets $B_n\in\goth C_n$ for each $n\in\Bbb N$.
Then the following holds.
\roster
\item"---" The restricted product $\eta^{\boldkey B}$ is equivalent to $\bigotimes_{n=1}^\infty\alpha_n$ if and only if 
$\eta_n\sim\alpha_n$ for each $n\in\Bbb N$ and
 $\sum_{n=1}^\infty H^2\big(\frac 1{\eta_n(B_n)}(\eta_n\restriction B_n),\alpha_n\big)<\infty$ \cite{Hi, Theorems~3.6}.
 \item"---" If $\beta$ be an MH-product of $(\eta_n)_{n=1}^\infty$ then there is $a>0$
 and a sequence $\boldkey C=(C_n)_{n=1}^\infty$ of Borel subsets such that
 $\beta= a\eta^{\boldkey C}$ \cite{Hi, Theorems~3.9}.
\endroster
\endproclaim

We also note that if $\eta^{\boldkey B}\sim\bigotimes_{n=1}^\infty\alpha_n$ then $\big(\bigotimes_{n=1}^\infty\alpha_n\big)(\boldsymbol B^k)>0$ for each $k>0$.
In particular, 
 $\prod_{n=1}^\infty\alpha_n(B_n)>0$.

\proclaim{Fact 2.5 \cite{Da, Proposition 1.6}} Let $S_n$ be a measure preserving invertible transformation of a standard probability space $(Y_n,\goth C_n,\eta_n)$, $n\in\Bbb N$.
Let $\boldkey B=(B_n)_{n=1}^\infty$ be a sequence of subsets $B_n\in\goth C_n$ with $0<\eta_n(B_n)<\infty$ for each $n\in\Bbb N$.
If  $\sum_{n=1}^\infty\frac{\eta_n(B_n\triangle S_nB_n)}{\eta_n(B_n)}<\infty$ then the infinite product transformation  $S:=\bigotimes_{n=1}^\infty S_n$ preserves $\eta^{\boldkey B}$.
\endproclaim

For completeness of our argument we reproduce here the proof of this fact from~\cite{Da}.

\demo{Proof} 
Fix $n>0$.
Take a Borel subset $A\subset Y_1\times\dots\times Y_n$ and put $A':=A\times B_{n+1}\times B_{n+2}\times\cdots$.
Then
$$
\eta^{\boldkey B}(SA')=\lim_{m\to\infty}(\eta^{\boldkey B}\restriction\boldkey B^m)(SA')=
\frac{\big(\bigotimes_{j=1}^n\eta_j\big)(A)}{\prod_{j=1}^n\eta(B_j)}\lim_{m\to\infty}\prod_{m>n}\frac{\eta_m(B_m\cap S_mB_m)}{\eta_m(B_m)}.
$$
It follows that $\eta^{\boldkey B}(SA')=\frac{(\bigotimes_{j=1}^n\eta_j)(A)}{\prod_{j=1}^n\eta(B_j)}=\eta^{\boldkey B}(A')$.
Hence  $\eta^{\boldkey B}\circ S=\eta^{\boldkey B}$, as desired.
\qed
\enddemo

We now state the main result of this section.
For that we will use the notation~ $\mu$, $(\mu_n)_{n=1}^\infty$ and $c(.)$  introduced in \S1.

\proclaim{Theorem 2.6} For each $n\in\Bbb N$, let
$a_n:=\frac 1{2n\log (n+1)}$ and let $\lambda_n\in(0,1)$ be defined by the formula
$c(\lambda_n)=\log (n+1)$.
Then $\Gamma\subset\text{\rom{Aut}}_1(X,\mu)$ and
the Poisson suspension $(X^*,\mu^*,\Gamma_*)$ of $(X,\mu,\Gamma)$ is sharply weak mixing of type $II_\infty$.
\endproclaim

\demo{Proof} It is easy to see that the function $c$ decreases on the interval $(0,1)$.
It follows that  $\lambda_n$ is well defined
and $\lambda_n\sim\frac 1{\sqrt{\log n}}$ as $n\to\infty$.
Fix $k\in\Bbb N$.
In view of \thetag{1-4} and as $|1-\lambda_n|<1$ and $a_n<1$,
$$
\sum_{n=1}^\infty \frac{\#(\gamma_k F_n\triangle F_n)}{\# F_n}|1-\lambda_n|a_n<\infty.
$$
Hence by Lemma~1.6(i), $T(\gamma_k)\in\text{Aut}_1(X,\mu)$ and $\chi(T(\gamma_k))=0$.
In a similar way, one can check that 
$ \big(\frac {d\mu}{d\mu\circ T(\gamma_k)^{-1}}\big)^2-1\in L^1(X,\mu)$ via Lemma~1.6(ii).
Furthermore, we deduce from  Lem\-ma~1.6(ii) that
$$
\align
\int_X\bigg(\bigg(\frac {d\mu}{d\mu\circ T(\gamma_k)^{-1}}\bigg)^2-1\bigg) d\mu
&\le 2\sum_{n=1}^k c(\lambda_n)a_n+
\sum_{n=k+1}^\infty\frac{c(\lambda_n)a_n}{n}\\
& =\sum_{n=1}^k \frac1n+
\sum_{n=k+1}^\infty\frac1{2n^2}\\
&=\log k + d +\overline o(1)
\endalign
$$
 as $k\to\infty$ for some $d\ge 0$.
This asymptotic inequality plus Lemma~1.5 imply that  $T_*$ is conservative
\footnote{Indeed, put $b_{\gamma_k}:=\frac 1{k\log k}$ for each $k\in\Bbb N$ in the statement of Lemma~1.5.}. Therefore Corollary~1.7  yields that $(X^*,\mu^*,T_*)$ is sharply weak mixing.

It remains to show that this system is of type $II_\infty$.
We recall that
 $(X^*,\nu^*, T_*)$ is the nonsingular direct product $\bigotimes_{n=1}^\infty (Z^*,\nu_n^*, L_*)$ and
$L_*$ preserves $\nu_n^*$ for each $n\in\Bbb N$.
According to \thetag{1-2}, for each $n\in\Bbb N$,
$$
\frac{d\nu_n^*}{d\mu_n^*}(\omega)=e^{\mu_n(A_n)-\nu_n(A_n)}
$$
at a.e.  $\omega\in[A_n]_0$.
Therefore, we obtain that for each $n\in \Bbb N$,
$$
H^2\bigg(\frac{1}{\nu_n^*([A_n]_0)}(\nu_n^*\restriction [A_n]_0),\mu_n^*\bigg)=
1-\frac{\sqrt{e^{\mu_n(A_n)-\nu_n(A_n)}}}{\sqrt{\nu_n^*([A_n]_0)}}\int_{[A_n]_0}d\mu_n^*=1-e^{-\frac{\mu_n(A_n)}{2}}.
$$
Since $\lambda_n\sim\frac1{\sqrt{\log n}}$ and $\sum_{n=1}^\infty\mu_n(A_n)=\sum_{n=1}^\infty\lambda_n\nu_n(A_n)=\sum_{n=1}^\infty\frac {\lambda_n}{2n\log( n+1)}<+\infty$, it follows that
$
\sum_{n=1}^\infty H^2\Big(\frac{1}{\nu_n^*([A_n]_0)}(\nu_n^*\restriction [A_n]_0),\mu_n^*\Big)<\infty.
$
Therefore Fact~2.4 yields that $\mu^*=\bigotimes_{n=1}^\infty\mu_n^*\sim \nu^{*\boldsymbol B}$,
where $\nu^{*\boldsymbol B}$ denotes  the restricted product of $(\nu_n^*)_{n=1}^\infty$ with respect
to ${\boldsymbol B}=([A_n]_0)_{n=1}^\infty$.
For each $\gamma\in\Gamma$, we have that $[A_n]_0\cap L(\gamma)_*[A_n]_0=[A_n]_0\cap [L(\gamma) A_n]_0=[A_n\cup L(\gamma)A_n]_0$
and hence
$$
\align
\sum_{n=1}^\infty\frac{\nu_n^*([A_n]_0\triangle L(\gamma)_*[A_n]_0)}{\nu_n^*([A_n]_0)}
&=
2\sum_{n=1}^\infty\bigg(1-\frac{\nu_n^*([A_n\cap L(\gamma) A_n]_0)}{\nu_n^*([A_n]_0)}\bigg)\\
&=
2\sum_{n=1}^\infty\Big(1-e^{-\frac12\nu_n(A_n\triangle L(\gamma) A_n)}\Big)<\infty
\endalign
$$
in view of \thetag{1-4}.
Therefore, by Fact~2.5, 
 $\nu^{*\boldsymbol B}\circ \gamma_*=\nu^{*\boldsymbol B}$.
Since 
$$
\prod_{n=1}^\infty \nu_n^*([A_n]_0)=\prod_{n=1}^\infty e^{-\nu(A_n)}=e^{-\sum_{n=1}^\infty\frac 1{2n\log (n+1)}}=0,
$$
we deduce from  Corollary~2.3 that $\nu^{*\boldsymbol B}(X^*)=\infty$. 
Thus $\nu^{*\boldsymbol B}$ is an infinite $\sigma$-finite $T_*$-invariant measure equivalent to $\mu^*$.
Hence $(X^*,\mu^*,T_*)$ is of type $II_\infty$.
\qed
\enddemo

\head 3. Poisson suspensions of type $III_0$
\endhead
Our purpose in this section is to prove the following theorem (as above, we use here the notation from \S1).

\proclaim{Theorem 3.1} Let  $(l_n)_{n=1}^\infty$ be a sequence of positive integers  such that $l_n\to\infty$, $l_1|l_2$, $l_2|l_3$,\dots and 
$\sum_{n=1}^\infty\frac{1}{n4^{l_n}}=+\infty$.
For each $n\in\Bbb N$, we let $\lambda_n:=2^{l_n}$ and $a_n:=\frac 1{2nc(\lambda_n)}$.
Then $T(\gamma)\in \text{\rom{Aut}}_1(X,\mu)$ for each $\gamma\in\Gamma$ and
the Poisson suspension $(X^*,\mu^*,T_*)$ of $(X,\mu, T)$ is sharply weak mixing of type $III_0$.
\endproclaim

\demo{Proof} We proceed in several steps.

{\sl Step 1.} We prove here  that $T(\gamma)\in \text{\rom{Aut}}_1(X,\mu)$  for each $\gamma\in\Gamma$ and $T_*$ is sharply weak mixing.
Fix $k\in\Bbb N$.
In view of \thetag{1-4} and as $1<\lambda_n\le c(\lambda_n)$ for each $n\in\Bbb N$,
$$
\sum_{n=1}^\infty \frac{\#(\gamma_k F_n\triangle F_n)}{\# F_n}|1-\lambda_n| a_n<\infty.
$$
Hence $T({\gamma_k})\in\text{Aut}_1(X,\mu)$ and $\chi(T({\gamma_k}))=0$ by Lemma~1.6(i).
In a similar way, one can check that 
$ \big(\frac {d\mu}{d\mu\circ  T({\gamma_k})^{-1}}\big)^2-1\in L^1(X,\mu)$ via Lemma~1.6(ii).
Furthermore, we deduce from  Lem\-ma~1.6(ii), the condition of the theorem  and \thetag{1-4} that
$$
\align
\int_X\bigg(\bigg(\frac {d\mu}{d\mu\circ T({\gamma_k})^{-1}}\bigg)^2-1\bigg) d\mu
&=\sum_{n=1}^\infty\frac 1{4n} \frac{\#(\gamma_k F_n\triangle F_n)}{\# F_n}
\\
& \le \sum_{n=1}^k \frac1{2n}+
\sum_{n=k+1}^\infty\frac1{4n^3}\\
&=\frac12\log k + d +\overline o(1)
\endalign
$$
as $k\to\infty$ for some $d\ge 0$.
It follows from this asymptotic inequality,  Lemma~1.5 and Corollary~1.7 that  
 $(X^*,\mu^*,T_*)$ is sharply weak mixing.

{\sl Step 2.} We show here  that $(X^*,\mu^*,T_*)$  is of type $III$.
Suppose, by contraposition, that $T_*$ is of type $II$.
Then there exists an ergodic $T_*$-invariant  $\mu^*$-equivalent $\sigma$-finite Borel measure $\vartheta$ on $X_*$.
We are going to show that $\vartheta$ is an MH-product of  $(\nu_n^*)_{n=1}^\infty$.
This can be deduced from \cite{Hi, Theorem~3.20}.
However we prefer an alternative---more direct---way, bypassing \cite{Hi, Theorem~3.20} whose proof is rather involved.
Fix $n\in\Bbb N$.
Since $\vartheta\sim\mu^*$ and $ \mu_1^*\otimes\cdots\otimes\mu_n^*\sim \nu_1^*\otimes\cdots\otimes\nu_n^*$,
the projection of  $\vartheta$ to $(Z^*)^{n}$ along the mapping
$$
\pi_n:X^*\ni (x_k^*)_{k=1}^\infty\mapsto (x_1^*,\dots,x_n^*)\in (Z^*)^n
$$
has the same collection of 0-measure  subsets as $\nu_1^*\otimes\cdots\otimes\nu_n^*$.
(We note that this projection is not, in general, $\sigma$-finite.) 
Therefore the disintegration of $\vartheta$ over  $\nu_1^*\otimes\cdots\otimes\nu_n^*$ along  $\pi_n$ is well defined (see, for instance, \cite{ChPo, Theorem~1}):
$$
\vartheta=\int_{(Z^*)^n}\delta_{y}\otimes\vartheta_y\,d(\nu_1^*\otimes\cdots\otimes\nu_n^*)(y),\tag3-1
$$
where $(Z^*)^n\ni y\mapsto\vartheta_y$ is the corresponding measurable field of $\sigma$-finite measures on the infinite product space $Z^*\times Z^*\times\cdots$.

We recall  that $R$ is a  totally dissipative free $\Gamma$-action on $Z$ that preserves $\nu_n$ for each $n\in\Bbb N$ (see \thetag{1-3}).
Hence for each $n$, the Poisson suspension $R_*:=(R(\gamma)_*)_{\gamma\in\Gamma}$ is a well defined weakly mixing 
 (it is Bernoulli)  $\nu_n^*$-preserving action of $\Gamma$.
 Moreover, for each $n\in\Bbb N$, 
 the actions $R_*$ and $L_*$ on $(Z^*,\nu_n^*)$ commute.
This action commutes with the Poisson suspension of $\Gamma$ on $(Z^*,\nu_n^*)$.
Hence,
 the  action $\Big(\underbrace{R(\gamma)_{*}\times \cdots \times R(\gamma)_{*}}_{n\text{ times}}\times I\Big)_{\gamma\in \Gamma}$ on $(X^*,\mu^*)$ commutes with $T_*$.
  Therefore, for each $n\in\Bbb N$ and $\gamma\in\Gamma$, there is $a_{n,\gamma}>0$
such that $\vartheta\circ (R(\gamma)_{*}\times \cdots \times R(\gamma)_{*}\times I)=a_{n,\gamma}\vartheta$.
Since  the transformation $R(\gamma)_{*}\times \cdots \times R(\gamma)_{*}\times I$ is conservative, $a_{n,\gamma}=1$.
On the other hand, it follows from \thetag{3-1} and the fact that $\nu^*_k\circ R(\gamma)_*=\nu_k^*$
for  $k=1,\dots,n$   that 
$$
\vartheta\circ(R(\gamma)_{*}\times \cdots \times R(\gamma)_{*}\times I)
=\int_{(Z^*)^n}\delta_{y}\otimes\vartheta_{R(\gamma)_{*}\times \cdots \times R(\gamma)_{*}y}\,d(\nu_1^*\otimes\cdots\otimes\nu_n^*)(y)
$$
each $\gamma\in\Gamma$.
Comparing this with \thetag{3-1} and using the uniqueness of  disintegration of $\vartheta$, we deduce that $\vartheta_y=\vartheta_{R(\gamma)_{*}\times \cdots \times R(\gamma)_*y}$ for 
$(\nu_1^*\otimes\cdots\otimes\nu_n^*)$-a.e. $y\in (Z^*)^n$.
Since the product $\Gamma$-action  $(R(\gamma)_*^{\times n})_{\gamma\in\Gamma}$ on $((Z^*)^n,\nu_1^*\otimes\cdots\otimes\nu_n^*)$ is ergodic, there is a $\sigma$-finite measure $\vartheta^{n}$ on $Z^*\times Z^*\times\cdots$ such that $\vartheta_y=\vartheta^n$ for a.e. $y\in(Z^*)^n$.
Then~\thetag{3-1} yields that $\vartheta=\nu_1^*\otimes\cdots\otimes\nu_n^*\times\vartheta^n$.
Since $n$ is arbitrary, we obtain that $\vartheta$ is a MH-product of  $(\nu_n^*)_{n=1}^\infty$.
Hence, by Fact~2.4,
there exists a sequence $\boldkey B=(B_n)_{n=1}^\infty$ of Borel subsets $B_n\subset X_n^*$ such that
$\vartheta$ is proportional to the restricted product $\nu^{*\boldkey B}$ of  $(\nu_n^*)_{n=1}^\infty$ with respect to ${\boldkey B}$.

If $\prod_{n=1}^\infty\nu_n^*(B_n)>0$ then $\nu^{*\boldkey B}\sim\bigotimes_{n=1}^\infty\nu_n^*$ by Corollary~2.3.
Hence
$\mu^*\sim\nu^*$. 
This implies that
 $\sqrt{\frac{d\mu}{d\nu}}-1\in L^2(X,\nu)$  (see \S1).
 However
 $$
 \int_X\bigg(\sqrt{\frac{d\mu}{d\nu}}-1\bigg)^2d\nu
  =\sum_{n=1}^\infty \int_Z\bigg(\sqrt{\frac{d\mu_n}{d\nu_n}}-1\bigg)^2d\nu_n
 =\sum_{n=1}^\infty(\sqrt{\lambda_n}-1)^2a_n=+\infty
 $$
because $a_n=\frac{1}{2nc(\lambda_n)}$, $c(\lambda_n)\sim\lambda_n^3$ as $n\to\infty$ and hence $\sum_{n=1}^\infty\frac{\lambda_n}{nc(\lambda_n)}=+\infty$ in view of the condition $\sum_{n=1}^\infty\frac1{n4^{l_n}}=+\infty$.
This contradiction implies that 
 $\prod_{n=1}^\infty\nu_n^*(B_n)=0$ or, equivalently,
$\sum_{n=1}^\infty\nu_n^*(Z^*\setminus B_n)=+\infty$.
On the other hand, it follows from our remark just below Fact~2.4 that $\prod_{n=1}^\infty\mu_n^*(B_n)>0$, i.e. $\sum_{n=1}^\infty\mu_n^*(Z^*\setminus B_n)<+\infty$.
Since $\lambda_n\ge 1$ and $\mu_n(A_n)\le 1$, we deduce from~\thetag{1-2} that
$\frac{d\nu_n^*}{d\mu_n^*}(\omega)\le e$ for $\nu_n^*$-a.e. $\omega\in Z^*$.
Hence 
$$
+\infty=\sum_{n=1}^\infty\nu_n^*(Z^*\setminus B_n)\le e^{}\sum_{n=1}^\infty\mu_n^*(Z^*\setminus B_n)<+\infty,
$$ 
a contradiction.
Hence $T_*$ is of type $III$.

{\sl Step 3.} We show here that $T_*$ is of type $III_0$.
Fix $n\in\Bbb N$.
According to~\thetag{1-1}, for each $\gamma\in\Gamma$ and $\mu^*$-a.e. $\omega=(\omega_k)_{k=1}^\infty\in X^*$,
$$
\frac{d\mu^*\circ T(\gamma)_*}{d\mu^*}(\omega)=\prod_{k=1}^\infty\frac{d\mu_k^*\circ L(\gamma)_*}{d\mu_k^*}(\omega_k).
$$
By \thetag{1-2} and the definition of $\mu_k$, 
$$
\align
\frac{d\mu_k^*\circ L(\gamma)_*}{d\mu_k^*}(\omega_k)&=\prod_{\omega_k(\{x_k\})=1}\frac{d\mu_k\circ L(\gamma)}{d\mu_k}(x_k)\\
&=
\prod_{\omega_k(\{x_k\})=1}\frac{f( \gamma x_k)}{f(x_k)}\in\{\lambda_k^n\mid n\in\Bbb Z\}
\endalign
$$
at $\mu_k^*$-a.e. $\omega_k\in Z^*$.
Moreover,  for a.e. $\omega_k\in[A_k]_0$, if $\gamma_*^{-1}\omega_k\in[A_k]_0$ for some $\gamma\in\Gamma$ then
$$
\frac{d\mu_k^*\circ L(\gamma)_*}{d\mu_k^*}(\omega_k)=
\prod_{\omega_k(\{x_k\})=1}\frac{f( \gamma x_k)}{f(x_k)}=1.
$$
We now set $B^n:=[A_1]_0\times\cdots\times[A_n]_0\times X_{n+1}^*\times X_{n+2}^*\times\cdots\subset X^*$.
Of course, $\mu^*(B^n)>0$.
It follows that each $\gamma\in\Gamma$ and  a.e. $\omega=(\omega_k)_{k=1}^\infty\in B^n$, if
$T(\gamma)_*\omega\in B^n$ then
$$
\frac{d\mu^*\circ T(\gamma)_*}{d\mu^*}(\omega)=\prod_{k=n+1}^\infty\frac{d\mu_k^*\circ L(\gamma)_*}{d\mu_k^*}(\omega_k)\in \{\lambda_{n+1}^m\mid m\in\Bbb Z\},\tag3-4
$$
because the multiplicative subgroup of $\Bbb R_*$ generated by $\lambda_k$ is a subgroup of the the multiplicative subgroup of $\Bbb R_*$ generated by $\lambda_{n+1}$
for each $k>n+1$.
It follows from~\thetag{3-4} and the definition of an essential value for the Radon-Nikodym cocycle of $T_*$ that $r(T_*)\subset \{\lambda_{n+1}^m\mid m\in\Bbb Z\}$.
Hence $r(T_*)\subset\bigcap_{n=1}^\infty \{2^{l_{n+1}m}\mid m\in\Bbb Z\}=\{1\}$.
Thus, $T_*$ is either of type $II$ or of type $III_0$.
The former contradicts  the assertion  proved on Step~2.
Hence $T_*$ is of type $III_0$.
\qed
\enddemo

\head 4. Poisson suspensions of type $III_\lambda$ with $0<\lambda<1$
\endhead
Fix $\lambda\in(0,1)$.
Our purpose in this section is to prove the following theorem (as above, we use here the notation from \S1).

\proclaim{Theorem 4.1} 
For $n\in\Bbb N$, let $\lambda_{2n-1}:=\lambda_{2n}^{-1}:=\lambda$ and 
$a_{2n-1}:=\lambda^{-1}a_{2n}:=\frac 1{n\log (n+1)}$.
 Then $T(\gamma)\in\text{\rom{Aut}}_1(X,\mu)$ for each $\gamma\in\Gamma$ and
the Poisson suspension $(X^*,\mu^*,T_*)$ of $(X,\mu,T)$ is sharply weak mixing of Krieger's type $III_\lambda$.
\endproclaim

\demo{Proof} Fix $k\in\Bbb N$.
As in the proof of Theorems~2.6 and 3.1,  it is straightforward to verify  that the series from Lemma~1.6(i) (with $\gamma_k$ in place of $\gamma$) converges.
Hence $T(\gamma_k)\in\text{Aut}_1(X,\mu)$ and $\chi(T(\gamma_k))=0$ by Lemma~1.6(i).
Since $c(\lambda^{-1})=\frac{c(\lambda)}\lambda$, it follows from the condition of the theorem that 
$$
c(\lambda_{2n})\nu(A_{2n})=c(\lambda_{2n-1})\nu(A_{2n-1})=\frac{c(\lambda)}{n\log(n+1)}\quad\text{ for all $n\in\Bbb N.$}
$$
Therefore Lemma~1.6(ii)  and \thetag{1-4} yield that
$$
\align
\int_X\bigg(\bigg(\frac {d\mu}{d\mu\circ T(\gamma_k)^{-1}}\bigg)^2-1\bigg) d\mu
&\le
\sum_{n=1}^{k} 2c(\lambda_n)a_n
+
\sum_{n=k+1}^{\infty} \frac{c(\lambda_n)a_n}n\\
&\le \sum_{n=1}^{\langle (k+1)/2\rangle}
\frac{2c(\lambda)}{n\log (n+1)}+
\sum_{n=\langle k/2\rangle}^\infty\frac{2c(\lambda)}{n^2\log(n+1)}\\
&=2c(\lambda)\log \log k + d +\overline o(1)
\endalign
$$
as $k\to\infty$ for some $d\ge0$.
Here, we use the notation $\langle .\rangle$ for the integer part.
This asymptotic inequality plus Lemma~1.5 and Corollary~1.7 imply that  $T_*$ is sharply weak mixing.\footnote{Put $b(\gamma_k):=\frac 1k$ in the statement of Lemma~1.5.}

Denote by $\Psi$ the action of the group $\bigoplus_{n=1}^\infty \Gamma$ on  $(X^*,\mu^*)$ generated by the transformations $\Psi(\gamma_1,\dots,\gamma_n):=L(\gamma_1)_*\times\cdots\times L(\gamma_n)_*\times I$, $\gamma_1,\dots,\gamma_n\in\Gamma$, $n\in\Bbb N$.
The dynamical system $(Z^*,\nu_n^*, L_*)$ is mildly mixing (it is  Bernoulli) for each $n\in\Bbb N$.
Hence, in view of Proposition~1.3, it suffices to prove that $\Psi$ is of type $III_\lambda$.

We first show that $\lambda^{-1}$ is an essential value of  the Radon-Nikodym cocycle 
of $\Psi$.
Fix $n>0$.
Since  $\sum_{k=n}^{+\infty}\mu_{2k+1}(A_{2k+1})=+\infty$ and $\mu_{2k+1}(A_{2k+1})\to 0$, there is $m_n>n$ such that 
$
\alpha_n:=\sum_{k=n}^{m_n-1}\mu_{2k+1}(A_{2k+1})\in (0.5,1).
$
Denote the probability space $\bigotimes_{k=2n+1}^{2m_n}(Z^*, \mu^*_k)$   
by $(Z^*_{2n+1,2m_n},\mu^*_{2n+1,2m_n})$.
We also let $\nu^*_{2n+1,2m_n}:=\bigotimes_{k=2n+1}^{2m_n}\nu^*_k$.
Since 
$$
\mu_{2n-1}(A_{2n-1})-\nu_{2n-1}(A_{2n-1})+\mu_{2n}(A_{2n})-\nu_{2n}(A_{2n})=0,
$$
it follows from \thetag{1-5} that for a.e. $(\omega_{2n-1},\omega_{2n})\in Z^*\times Z^*$,
$$
\aligned
\frac{d(\mu_{2n-1}^*\times \mu_{2n}^*)}{d(\nu_{2n-1}^*\times\nu_{2n}^*)}(\omega_{2n-1},\omega_{2n})
&=\sum_{k=0}^\infty\lambda^k1_{[A_{2n-1}]_k}(\omega_{2n-1})
\sum_{k=0}^\infty\lambda^{-k}1_{[A_{2n}]_k}(\omega_{2n})\\
&=\sum_{k=-\infty}^\infty\lambda^k1_{B_{n,k}}(\omega_{2n-1},\omega_{2n}),
\endaligned
\tag4-1
$$
where $B_{n,k}=\sum_{k=j-r} 1_{[A_{2n-1}]_j}1_{[A_{2n}]_r}$.
Hence the mapping
$$
\vartheta_n:Z^*_{2n+1,2m_n}\ni \omega\mapsto\log_\lambda\frac{\mu^*_{2n+1,2m_n}}{\nu^*_{2n+1,2m_n}}(\omega)\in\Bbb Z
$$
is well defined.
Applying \thetag{4-1}, we obtain that
$$
\vartheta_n(\omega)=\sum_{k=n+1}^{m_n}(\omega(A_{2k-1})-\omega(A_{2k})),\quad\omega\in Z^*_{2n+1,2m_n}.
$$
Thus, $\mu^*_{2n+1,2m_n}\circ\vartheta_n^{-1}$ is the distribution of the difference of two independent Poisson random variables, one with parameter $\alpha_n$, the other with
parameter $\lambda\alpha_n$.
In other words, $\mu^*_{2n+1,2m_n}\circ\vartheta_n^{-1}$ is the 
Skellam distribution with parameters $\alpha_n,\lambda\alpha_n$.
For $i=0,1$, let $\Delta_i:=\vartheta_n^{-1}(\{i\})\subset Z^*_{2n+1,2m_n}$.
Then \footnote{Skellam (1946) and Prekopa (1953) represented the Skellam distribution using the modified Bessel function of the first kind. The result we are referring to is a direct consequence of this and standard facts on Bessel functions \cite{AbSt, pp. 374--378}.}
$$
\mu^*_{2n+1,2m_n}(\Delta_1)=
e^{-\alpha_n(1+\lambda)}\sum_{k=0}^\infty\frac{\alpha_n^{k+1}(\lambda\alpha_n)^k}{(k+1)!k!}>\frac{\alpha_n}{e^{\alpha_n(1+\lambda)}}>\frac1{16}
    $$
    and $\mu^*_{2n+1,2m_n}(\Delta_0)>\mu^*_{2n+1,2m_n}(\Delta_1)$.
Therefore,
$$
\frac{\nu^*_{2n+1,2m_n}(\Delta_0)}{\nu^*_{2n+1,2m_n}(\Delta_1)}=
\frac{\lambda\mu^*_{2n+1,2m_n}(\Delta_0)}{\mu^*_{2n+1,2m_n}(\Delta_1)}\ge \lambda.
$$
Hence, we can select a subset $\Delta'_1\subset\Delta_1$  and a transformation
$S$ from the full group of the action $\bigoplus_{j=2n+1}^{2m_n} L_*$ defined on the space $(Z^*_{2n+1,2m_n},
\mu^*_{2n+1,2m_n})$
such that $\nu^*_{2n+1,2m_n}(\Delta'_1)=\lambda\nu^*_{2n+1,2m_n}(\Delta_1)$ and  $S\Delta'_1\subset\Delta_{0}$.
Then
$$
\mu^*_{2n+1,2m_n}(\Delta_1')=\lambda\nu^*_{2n+1,2m_n}(\Delta_1')=\lambda^2\nu^*_{2n+1,2m_n}(\Delta_1)=\lambda\mu^*_{2n+1,2m_n}(\Delta_1)>\frac\lambda{16}.
$$
If $\omega\in \Delta_1'$ then  $S\omega\in \Delta_0$ and
$$ 
\frac{d\mu^*_{2n+1,2m_n}\circ S}{d\mu^*_{2n+1,2m_n}}(\omega)=
\frac{d\mu^*_{2n+1,2m_n}}{d\nu^*_{2n+1,2m_n}}(S\omega)\frac{d\nu^*_{2n+1,2m_n}}{d\mu^*_{2n+1,2m_n}}(\omega)=\lambda^{-1}.
$$
For $l\in\Bbb N$ and a Borel subset $C\subset (Z^*)^l$, we denote by $[C]_l$ the cylinder $C\times Z^*\times Z^*\times \cdots\subset X^*$.
Now, for a cylinder    $[C]_{2n}$ in $X^*$, we have that
\roster
\item"---"
$[C\times\Delta'_1]_{2m_n}\subset [C]_{2n}$,
\item"---"
$\mu^*([C\times\Delta'_1]_{2m_n})=\mu^*([C]_{2n})\mu^*_{2n+1,2m_n}(\Delta_1')>\frac\lambda {16}\mu^*([C]_{2n})$,
\item"---"
the transformation $I\times S\times I$ belongs to the full group  of $\Psi$ and 
$$
\frac{d\mu^*\circ (I\times S\times I)}{d\mu^*}(\omega)=\lambda^{-1}\qquad\text{for a.e. $\omega\in[C\times\Delta'_1]_{2m_n}$.}
$$ 
\endroster
Since the family of cylinders $\{[C]_{2n}\mid C\subset (Z^*)^{2n}, n\in\Bbb N\}$ is dense
in the entire Borel $\sigma$-algebra on $X^*$, it follows  that $\lambda^{-1}\in r(\Psi)$.
Since the Radon-Nikodym cocycle of the system $(X^*,\mu^*,\Psi)$ takes its values in the subgroup $\{\lambda^n\mid n\in\Bbb Z\}$, it follows that
$\Psi$ is of type $III_\lambda$, as desired.
\qed
\enddemo

\head 5. Poisson suspensions of type $III_1$ 
\endhead
We note that some examples  type $III_1$ ergodic Poisson suspensions 
over a totally dissipative $\Bbb Z$-actions were  given in \cite{DaKoRo2}.
Now we present an alternative construction which is a natural modification of  the $III_\lambda$-construction from \S4.
It works well for arbitrary countable amenable groups.

Fix $\lambda_1,\lambda_2\in(0,1)$ such that $\log\lambda_1$ and $\log\lambda_2$ are rationally independent.
 The following theorem  follows from Theorem~4.1 (as above, we use here the notation from \S1).

\proclaim{Theorem 5.1} 
For $n\in\Bbb N$, let $\lambda_{4n-3}:=\lambda_{4n-2}^{-1}:=\lambda_1$, 
$\lambda_{4n-1}:=\lambda_{4n}^{-1}:=\lambda_2$,
$a_{4n-3}:=\lambda^{-1}_1a_{4n-2}:=\frac 1{n\log (n+1)}$ 
and $a_{4n-1}:=\lambda_2^{-1}a_{4n}:=\frac 1{n\log (n+1)}$.
Then $T(\gamma)\in\text{\rom{Aut}}_1(X,\mu)$ for each $\gamma\in\Gamma$ and
the Poisson suspension $(X^*,\mu^*,T_*)$ of $(X,\mu,T)$ is sharply weak mixing of type $III_1$.
\endproclaim

\demo{Idea of the proof} The proof of the conservativeness of $T_*$ is only a slight obvious modification of the first part of the proof of Theorem~4.1.
To complete the proof, one show follow the proof of Theorem~4.1 and to use the fact that the direct product of two nonsingular $\bigoplus_{n=1}^\infty\Gamma$-actions $\Psi_1$ and $\Psi_2$, one is of type $III_{\lambda_1}$, the other one is of type $III_{\lambda_2}$, considered as an action of the group $(\bigoplus_{n=1}^\infty\Gamma)\times
(\bigoplus_{n=1}^\infty\Gamma)$,  is of type $III_1$.
\qed
\enddemo

\enddocument